\documentclass[12pt]{article}
\usepackage{amsmath}
\newcommand{\be}{\begin{equation}}
      \newcommand{\ee}{\end{equation}}
      \newcommand{\ba}{\begin{eqnarray}}
       \newcommand{\ea}{\end{eqnarray}}
\newcommand{\ban}{\begin{eqnarray*}}
       \newcommand{\ean}{\end{eqnarray*}}

 \newcommand{\qed}{\hspace*{\fill}\rule{3mm}{3mm}\quad}
 \newcommand{\Pf}{\noindent {\em Proof.} }

\newcommand{\sect}[1]{\section{#1} \setcounter{equation}{0}}

\newtheorem{lem}{Lemma}[section]

\begin{document}
 \newtheorem{defn}[lem]{Definition}
 \newtheorem{theo}[lem]{Theorem}
 \newtheorem{prop}[lem]{Proposition}
 \newtheorem{rk}[lem]{Remark}
 \newtheorem{ex}[lem]{Example}
 \newtheorem{note}[lem]{Note}
 \newtheorem{conj}[lem]{Conjecture}

\title{Curvature Estimates for the Ricci Flow I}
\author{Rugang Ye \\ {\small Department of Mathematics} \\
{\small University of California, Santa Barbara}}
\date{}
\maketitle

\sect{Introduction}

In this paper we present several curvature estimates for
solutions of the Ricci flow
\ba
\frac{\partial g}{\partial t}=-2Ric.
\ea
These estimates depend on the smallness of
certain local $L^{\frac{n}{2}}$ integrals of
the norm of the Riemann curvature tensor $|Rm|$, where 
$n$ denotes the dimension of the manifold.
A key common property of these integrals is
scaling invariance, thanks to the critical
exponent $\frac{n}{2}$.  Because of this property,
they are very natural
and contain particularly
rich geometric information. (Note that in dimension 4, the
norm square of the Riemann curvature tensor is closely related to the
Gauss-Bonnet-Chern integrand.)

To formulate our results, we need some terminologies.  Consider a  
Riemannian manifold $(M, g)$ ($g$ denotes the metric) possibly
with boundary.  For convenience, we define the distance between two points of $M$ to be $\infty$, if they belong to 
two different connected components.  Consider a point $x\in M$.
If $x$ is in the interior of $M$, we define the distance $d(x, \partial M)=d_g(x, \partial M)$ to be
$\sup\{r>0: B(x,r) \mbox{ is compact}$  and contained in the interior of  $M\}$, where
$B(x, r)$ denotes the closed geodesic ball of center $x$ and radius $r$. If $M$ has
a boundary and $x \in \partial M$, then $d(x, \partial M)$ is the ordinary distance
from $x$ to $\partial M$ and equals zero. (For example, $d(x, \partial M)=\infty$ if
$M$ is closed.)    \\

\noindent {\bf Notations} Let $g=g(t)$ be a family of metrics on $M$. Then $d(x, y, t)$ denotes the distance
between $x, y \in M$ with respect to the metric $g(t)$, and $B(x, r ,t)$ denotes the closed geodesic ball of center $x\in M$ and radius $r$ with
respect to the metric $g(t)$. The volume of $B(x, r, t)$ with respect to $g(t)$ will 
often be denoted by $V(x,r,t)$. We shall often use $dq$ to denote $dvol_{g(t)}$.
These notations naturally extend when $g$ and (or) $t$ are replaced by 
other notations.\\

We set $\alpha_n=\frac{1}{40(n-1)}$, $\epsilon_0=\frac{1}{42}$ and
$\epsilon_1=\frac{\epsilon_0}{\sqrt{1+2\alpha_n\epsilon_0^2}}$.
(These constants are
not meant to be optimal. One can improve them by carefully examing
the proofs.)  \\

Our first result involves the concept of $\kappa$-noncollapsedness due
to Perelman [P].\\

\noindent {\bf Definition} Let $g$ be a Riemannian metric on a manifold $M$ of dimension
$n$. Let $\kappa$ and $\rho$ be positive numbers. We say that $g$ is {\it $\kappa$-noncollapsed on the scale} 
$\rho$,
if $g$ satisfies $vol(B(x, r))\ge \kappa r^n $
for all $x \in M$ and $r>0$ with the properties $r < \rho$ and $\sup\{|Rm|(x): x\in B(x, r)\}
\le r^{-2}$.  We say that a family of Riemannian metrics $g=g(t)$ is $\kappa$-noncollapsed on the scale 
$\rho$, if $g(t)$ is $\kappa$-noncollapsed on the scale $\rho$ for each $t$ (in the given domain). 
\\

By [Theorem 4.1, P], a smooth solution $g$ of the Ricci flow on $M \times [0, T)$
for a closed manifold $M$ and a finite $T$
is $\kappa$-noncollapsed on the scale $\sqrt{T}$, where $\kappa$ depends on
the initial metric and (an upper bound of $T$).  For  a variant of [Theorem 4.1, P] involving the scalar curvature,  
which implies [Theorem 4.1, P], see Theorem A.1 in Appendix A. By Theorem A.1, $g$ is 
$\kappa$-noncollapsed on the scale $\rho$ for an arbitary positive number $\rho$, where 
$\kappa$ depends on the initial metric and (an upper bound of) $T+\rho^2$.  \\

\noindent {\bf Theorem A} {\it  For each  positive number $\kappa$ and each  natural number $n\ge 3$
there are positive constants
$\delta_0=\delta_0(\kappa, n)$ and $C_0=C_0(n,\kappa)$  depending only on $\kappa$ and
$n$ with the following property. 
Let $g=g(t)$ be a smooth solution of the Ricci flow on
$M \times [0, T)$ for a manifold $M$ of dimension $n\ge 3$ and some
(finite or infinite) $T>0$,
which is $\kappa$-noncollapsed  on
the scale $\rho$ for some $\kappa>0$ and $\rho>0$.
Consider $x_0 \in M$ and $0<r_0\le \rho$, which satisfy
$r_0<d_{g(t)}(x_0, \partial M)$ for
each $t \in [0, T)$. Assume that
\ba \label{smallRmintegralA-1}
\int_{B(x_0, r_0, t)} |Rm|^{\frac{n}{2}}(\cdot, t)dq(\cdot, t) \le \delta_0
\ea
for all $t \in [0, T)$.
Then we have
\ba \label{RmboundA-1}
|Rm|(x, t) \le \alpha_n t^{-1}+(\epsilon_0 r_0)^{-2}
\ea
whenever $t \in (0, T)$ and $d(x_0, x, t)<\epsilon_0 r_0$, and  
\ba \label{RmboundintA-1}
|Rm|(x, t) \le C_0 \max\{\frac{1}{r_0^{2}}, \frac{1}{t}\}\sup\limits_{0\le s \le t}\left(\int_{B(x,\frac{1}{2}r(t), s)} |Rm|^{\frac{n}{2}}(\cdot, s)dq(\cdot, s)\right)^{\frac{2}{n}}
\ea
whenever $0<t<T$ and $d(x_0, x, t)\le \frac{1}{2}r(t),$
where $r(t)=\epsilon_1\min\{r_0, \sqrt{t}\}.$ (Obviously, it follows  that if the assumptions hold
on $[0, T]$, then the estimates (\ref{RmboundA-1}) and
(\ref{RmboundintA-1}) hold on
$[0, T]$. This remark also applies to the results below.) } \\

Note that the constant $\delta_0$ depends on $n$ decreasingly and depends on $\kappa$ increasingly, i.e. 
$\delta_0(n, \kappa)$ is a decreasing function of $n$ and an increasing function of $\kappa$. In contrast,  the 
constant $C_0$ depends on $n$ increasingly and depends on $\kappa$ decreasingly. 
The dependences of the constants in Theorem B and Theorem C are of similar nature.\\

The following corollary is a consequence of Theorem A and Theorem A.1. \\

\noindent {\bf Corollary A} {\it Let $g=g(t)$ be a smooth solution of
the Ricci flow  on $M \times [0, T)$ for a closed
manifold $M$  of dimension $n\ge 3$  and some finite $T>0$. Let $\rho>0$ be 
a positive number.
 There is a positive constant $\delta_0=
\delta_0(n, T+\rho^2, g(0))$ depending
on $n, T+\rho^2$ and $g(0)$ with the following property.
Assume that (\ref{smallRmintegralA-1}) holds true
for all $x_0\in M$, all $t \in [0,T)$, and some $r_0$ satisfying
$r_0 \le \rho$ for
all $t \in [0, T)$.  Then
(\ref{RmboundA-1}) holds true for all $x \in M$ and $t \in [0, T)$.
Consequently,  $g(t)$ extends to a smooth solution
of the Ricci flow over $[0,T']$ for some $T'>T$. } \\

A corollary of Theorem A for the case $T=\infty$ analogous to Corollary B 
(under the assumption of $\kappa$-noncollapsedness) also holds true. We omit the statement. 
Our second result does not involve the condition of $\kappa$-noncollapsedness.
Instead, smallness of $L^{\frac{n}{2}}$
integrals of $Rm$
over balls of varying center and radius measured against a volume ratio
is assumed. \\

\noindent {\bf Theorem B} {\it
For each natural number $n\ge 3$ there are positive constants $
\delta_0=\delta_0(n)$ and $C_0=C_0(n)$ depending only on $n$ with the following property.
Let $g=g(t)$ be a smooth solution of the Ricci flow on
$M \times [0, T)$ for a manifold $M$ of dimension $n \ge 3$ and some
(finite or infinite) $T>0$.
Consider $x_0 \in M$ and $r_0>0$,which satisfy
$r_0\le diam_{g(t)}(M, x_0)$ and $r_0<d_{g(t)}(x_0, \partial M)$ for
each $t \in [0, T)$.
Assume that
\ba \label{smallRmintegralB-1}
\int_{B(x, r, t)} |Rm|^{\frac{n}{2}}(\cdot,t)dq(\cdot, t) \le \delta_0\frac{V(x, r, t)}{r^n}
\ea
whenever $t \in [0, T), 0<r\le \frac{r_0}{2}$ and $x \in B(x_0, \frac{r_0}{2}, t)$.  Then we have
\ba \label{RmboundB-1}
|Rm|(x, t) \le \alpha_n t^{-1}+(\epsilon_0 r_0)^{-2}
\ea
whenever $t \in (0, T)$ and $d(x_0, x, t)<\epsilon_0 r_0$, and 
\ba \label{RmboundintB-1}
|Rm|(x, t) \le C_0  \sup\limits_{0\le s \le t}\left(
\frac{\int_{B(x, \frac{1}{9}r(t), s)} |Rm|^{\frac{n}{2}}(\cdot, s)dq(\cdot, s)}{V(x,
\frac{1}{9}r(t), s)}\right)^{\frac{2}{n}},
\ea
whenever $0<t<T$ and $d(x_0, x, t)\le \frac{1}{2} \epsilon_1\min\{r_0, \sqrt{t}\}$, where
$r(t)=\epsilon_1 \min\{r_0, \sqrt{t}\}.$
} \\

\noindent {\bf Corollary B} {\it Let $g=g(t)$ be a smooth solution of
the Ricci flow on $M \times [0, T)$ for a 
manifold $M$ of dimension $n \ge 3$ and some (finite or infinite) $T>0$,such that $g(t)$ is complete
for each $t \in [0, T)$.  Assume that (\ref{smallRmintegralB-1}) holds true
for all $x_0\in M$, all $t \in [0,T)$, and all $0<r\le r_0$ for
some positive number $r_0$ satisfying $r_0\le diam_{g(t)}(M)$
for all $t \in [0, T)$. Then
(\ref{RmboundB-1}) and (\ref{RmboundintB-1})  hold for all $x \in M$ and $t \in [0, T)$.
Consequently,  $g(t)$ extends to a smooth solution
of the Ricci flow over $[0,T']$ for some $T'>T$ if $T$ is finite. If
$T=\infty$, then $g(t)$ subconverges smoothly as $t \rightarrow T$.} \\

Our third result  does not involve the condition of $\kappa$-noncollapsedness,
and employs only a fixed center and a fixed radius for $L^{\frac{n}{2}}$ integrals of the
norm of the Riemann curvature tensor. But a lower bound for the
Ricci curvature is assumed. \\

\noindent {\bf Theorem C} {\it
For each natural number $n\ge 3$ there are positive constants $
\delta_0=\delta_0(n)$ and $C_0=C_0(n)$  depending only on $n$ with the following property.
Let $g=g(t)$ be a smooth solution of the Ricci flow on
$M \times [0, T)$ for a manifold of dimension $n\ge 3$ and some
(finite or infinite) $T>0$.  Consider $x_0 \in M$ and $r>0$, which satisfy
$r_0\le diam_{g(t)}(M)$ and $r_0<d_{g(t)}(x_0, \partial M)$ for
each $t \in [0, T)$. Assume that
\ba \label{riccibound}
Ric(x, t) \ge -\frac{n-1}{r_0^2}g(x,t)
\ea
whenever $t \in [0, T)$ and $d(x_0, x, t)\le r_0$ ($g(x,t)=g(t)(x)$
and $Ric(x,t)$ is the Ricci tensor of $g(t)$ at $x$), and that
\ba \label{smallRmintegralC-1}
\int_{B(x_0, r_0, t)} |Rm|^{\frac{n}{2}}(\cdot,t) dq(\cdot, t)\le \delta_0\frac{V(x_0, r_0, t)}{r_0^n}
\ea
for all $t \in [0, T]$.  Then we have
\ba \label{RmboundC-1}
|Rm|(x, t) \le \alpha_n t^{-1}+(\frac{1}{2}\epsilon_0 r_0)^{-2}
\ea
whenever $t \in (0, T)$ and $d(x_0, x, t)<\frac{1}{2}\epsilon_0 r_0$,
and
\ba \label{RmboundintC-1}
|Rm|(x, t) \le C_0  \sup\limits_{0\le s  \le t}\left(
\frac{\int_{B(x, \frac{1}{9}r(t), s)} |Rm|^{\frac{n}{2}}(\cdot, s)dq(\cdot, s)}{V(x,
\frac{1}{9}r(t), s)}\right)^{\frac{2}{n}}
\ea
whenever $0<t<T$ and $d(x_0, x, t)\le \frac{1}{2}r(t)$, where $r(t)= \epsilon_1 \min\{\frac{1}{2} r_0, \sqrt{t}\}$.   
We also have 
\ba \label{RmboundintC-2}
|Rm|(x, t) \le C_0  \left(\frac{r_0}{r(t)}\right)^2 \sup\limits_{0\le s  \le t}\left(
\frac{\int_{B(x_0, r_0, s)} |Rm|^{\frac{n}{2}}(\cdot, s)dq(\cdot, s)}{V(x_0,
r_0, s)}\right)^{\frac{2}{n}}
\ea
whenever $0<t<T$ and $d(x_0, x, t)\le \frac{1}{2} r(t).$  }
\\

An  elliptic analogue of (\ref{RmboundintC-2})  for Einstein metrics can be
found in [An]. Obviously, Theorem C and its corresponding version for the 
modified  Ricci flow (see discussions below) can also be applied to 
Einstein metrics. 
 \\

\noindent {\bf Corollary C} \noindent {\it Let $g=g(t)$ be a smooth solution of
the Ricci flow on $M \times [0, T)$ for a 
manifold $M$ of dimension $n\ge 3$ and some (finite or infinite) $T>0$, such that $g(t)$ is complete
for each $t\in [0, T)$. Assume that (\ref{riccibound}) and (\ref{smallRmintegralC-1}) hold true
for all $x_0\in M$ and
some positive number $r_0$ satisfying $r_0\le rad_{g(t)}(M)$
for all $t \in [0, T)$. Then
(\ref{RmboundC-1}) and (\ref{RmboundintC-1}) holds for all $x \in M$ and $t \in [0, T)$.
Consequently, $g(t)$ extends to a smooth solution
of the Ricci flow over $[0,T']$ for some $T'>T$, if $T$ is finite. If
$T=\infty$, then $g(t)$ subconverges smoothly as $t \rightarrow T$.} \\

Note that the condition $r_0 \le diam_{g(t)}(M)$ appears in Theorem B and Theorem C, but not 
in Theorem A.  \\

The above results extend to  
the modified Ricci flow 
\ba \label{chi}
\frac{\partial g}{\partial t}=-2Ric+\lambda(g, t) g
\ea
with a scalar function $\lambda(g, t)$ independent of $x \in M$. 
The volume-normalized 
Ricci flow 
\ba \label{volumenormalize}
\frac{\partial g}{\partial t}=-2Ric+\frac{2}{n} {\hat R} g
\ea
on a closed manifold, with $\hat R$ denoting the average scalar curvature, is  an example of the modified Ricci flow.
We present two 
extentions of the curvature estimates to the modified Ricci flow.   \\

\noindent {\bf Extention I} {\it Theorem A, Theorem B and Theorem C hold true for the modified Ricci flow (\ref{chi}), with the modification that 
the constants $\delta_0$ and $C_0$ in each theorem 
depend in addition on $r_0^2 |\min\{\inf_{[0, T)} \lambda, 0\}|$ which is assumed to be finite.
(In other words, $\delta_0$ and $C_0$ depend in addition on a nonpositive lower bound of $r_0^2 \lambda$. )
 } \\

\noindent {\bf Extention II} {\it Theorem A, Theorem B and Theorem C in the case $T<\infty$ hold true for the 
modified Ricci flow (\ref{chi}), with the modification that the constants $\delta_0$ and $C_0$ in 
each theorem depend in addition on $|\min\{\inf_{0\le t_1 <t_2 <T} \int_{t_1}^{t_2} \lambda, 0\}|$ 
which is assumed to be finite.  (In other words, $\delta_0$ and $C_0$ depend in addition on a nonpositive 
lower bound of $\int_{t_1}^{t_2} \lambda$.) } \\ 

In both extentions, the dependence of $\delta_0$ is decreasing, and the dependence of $C_0$ is increasing.   Extention I can be proved by directly adapting the proofs of Theorem A, Theorem B and Theorem C. Extention II 
can be proved by converting the modified Ricci flow into the Ricci flow, applying Theorem A, Theorem B and 
Theorem C, and then converting the obtained estimates back to the modified Ricci flow.

 Similar results also hold for many other evolution  equations. This will be presented elsewhere.

The results in this paper were obtained some time ago.

Analogous results involving other types of $L^p$
 integrals of $|Rm|$, including the case
 $p<\frac{n}{2}$ and space-time integrals,
 will
be presented in sequels of this paper. In particular, the case of space-time integrals 
is presented in [Ye3].

\sect{A Linear Parabolic Estimate}

In this section we present a linear parabolic estimate based on
Moser's iteration, which will be needed for establishing our
curvature estimates.
First we fix some notations. Consider a Riemannian manifold $(M, g)$ of dimension
$n$, possibly with boundary. Let $\Omega$ be a domain in $M$.  The Sobolev constant $C_{S, g}(\Omega)$ is defined to be the smallest number $C_{S, g}(\Omega)$ such that
\ba \label{sobolev-1}
\|f\|_{\frac{n}{n-1}}
\le C_{S,g}(\Omega) \|\nabla f\|_1
\ea
for all Lipschitz functions $f$ on $\Omega$ with compact support contained in the interior of 
$\Omega$ (i.e. $\Omega-\partial M$), where
$\| \cdot \|_{p}$ means the $L^p$-norm. More precisely,
$$
C_{S, g}(\Omega)=\sup\{\|f\|_{\frac{n}{n-1}}: f \in C^1_c(\Omega), \|\nabla f\|_1=1\},
$$
where $C^1_c(\Omega)$ is the space of $C^1$ functions on $\Omega$ with compact support contained 
in the interior of $\Omega$. As is well-known,
$C_{S,g}(\Omega)$ equals the isoperimetric constant $C_{I,g}(\Omega)$,
which is defined to be
$
\sup\{\frac{vol(\Omega')^{\frac{n-1}{n}}}{vol(\partial \Omega')}:
\Omega' \subset \Omega \mbox{ is a } C^1 \mbox{ domain in } \Omega 
\mbox{ with compact closure}.\}
$
The $L^2$-Sobolev constant $C_{S,2,g}(\Omega)$ is defined to be the smallest
number $C_{S,2,g}(\Omega)$ such that
$$
\|f\|_{\frac{2n}{n-2}} \le C_{S,2,g}(\Omega)\|\nabla f\|_2.
$$
There holds 
\ba \label{sobolev-2}
C_{S,2,g}(\Omega)\le \frac{2(n-1)}{n-2} C_{S, g}(\Omega).
\ea
 Indeed, applying (\ref{sobolev-1}) to $|f|^{\frac{2(n-1)}{n-2}}$ we deduce
\ban
(\int_{\Omega} |f|^{\frac{2n}{n-2}}dvol_{g} )^{\frac{n-1}{n}} &\le& C_{S,g}(\Omega)\int|\nabla
|f|^{\frac{2(n-1)}{n-2}}| dvol_g  \\ &\le& C_{S,g}(\Omega) \frac{2(n-1)}{n-2} \|\nabla f\|_2
(\int_{\Omega} |f|^{\frac{2n}{n-2}}dvol_g)^{\frac{1}{2}}.
\ean
Hence the claimed inequality follows.

It is easy to see that $C_{S, g}, C_{S, 2, g}$ and $C_{I, g}$ are nondecreasing, i.e. for example
$C_{S,2, g}(\Omega_1) \le C_{S,2,g}(\Omega_2)$ if $\Omega_1 \subset \Omega_2$.  

We set $C_{S, g}(\bar \Omega)=C_{S, g}(\Omega)$, $C_{S,2,g}(\bar \Omega)=C_{S, 2 ,g}(\Omega)$
and $C_{I, g}(\bar \Omega)=C_{I, g}(\Omega)$.

The following result is taken from [Ye2]. We include
 the proof for the
convenience of the reader, and for the reason of
verifying the explicit dependence on
the Sobolev constant, which
is important for the
curvature estimates in this paper.

\begin{theo}
 Let $M$ be a smooth manifold of dimension $n$ and $g=g(t)$ a smooth family of
Riemannian metrics on $M$ for $t \in [0, T]$. Let $f$ be a nonnegative Lipschitz continuous function on $M \times [0,T]$ satisfying
\ba \label{linearequation1}
\frac{\partial f}{\partial t} \leq \Delta f + af
\ea
on $M \times [0, T]$ in the weak sense, where $a$ is a nonnegative constant
and $\Delta=\Delta_{g(t)}$. Let $x_0$ be an interior point of $M$. Then
we have for each $p_0>1$ and $0<R< d_{g(0)}(x_0, \partial M)$
\ba \label{linearestimate1}
|f(x,t)|&\leq& (1+\frac{2}{n})^{\frac{\sigma_n}{p_0}} C_{S,2}^{\frac{n}{p_0}}\left(
ap_0+\frac{\gamma}{2} +\frac{n}{2}(1+\frac{n}{2})^2 \cdot \frac{1}{t}+
\frac{(n+2)^2e^{-\lambda_* T}}{4R^2}\right) ^{n+2\over 2p_0} \nonumber \\
&\cdot&
\Bigl(\int^T_0\int\limits_{B(x_0, R,0)} f^{p_0}(\cdot, t)dvol_{g(t)}dt\Bigr)^{\frac{1}{p_0}}~,
\ea
whenever $0<t\le T$ and $d_{g(0)}(x_0, x) \le \frac{R}{2}$, where $\sigma_n=\sum_0^{\infty} \frac{2k}{(1+\frac{2}{n})^k}$, $\gamma$ denotes the maximum value of the trace of $\frac{\partial g}{\partial t}$ on
$B(x_0, R, 0)\times [0, T]$,
$\lambda_*$ denotes the minimum eigenvalue of $\frac{\partial g}{\partial t}$
on $B(x_0, R, 0) \times [0, T]$, and
$$C_{S,2}=\max\limits_{0\le t\le T} C_{S,g(t), 2}(B(x_0, R, 0)).$$

The same estimate holds if we replace $B(x_0, R, 0)$ by $B(x_0, R, T)$, and
$-\lambda_*$ by $\lambda^*$, which denotes the maximum eigenvalue of
$\frac{\partial g}{\partial t}$
on $B(x_0, R, T) \times [0, T]$.
\end{theo}
\Pf We handle the case of $B(x_0, R, 0)$, while the other case is similar.
Let $\eta$ be a non-negative Lipschitz function on $M$ whose
support is contained in $B(x_0, R, 0)$.
 The partial differential inequality
(\ref{linearequation1}) implies for $p\geq 2$
\[
\frac{1}{p}\frac{\partial}{\partial t}\int f^p\eta^2
dvol_{g(t)}\leq -\int\nabla (\eta^2f^{p-1}) \cdot \nabla f\,dvol_{g(t)}+\int
bf^p\eta^2dvol_{g(t)}+\frac{1}{p}\int f^p\eta^2\frac{\partial}{ \partial t}dvol_{g(t)}~.
\]
We'll omit the notation $dvol_{g(t)}$ below. We have
\ban
-\int \nabla(\eta^2f^{p-1}) \cdot \nabla f\, & = &
-\frac{4(p-1)}{p^2}\int |\nabla (\eta f^{p/2})|^2
+\frac{4}{ p^2}\int |\nabla \eta|^2 f^p   \\
& & +\frac{4(p-2)}{ p^2}\int \nabla (\eta f^{p/2})
f^{p/2}\nabla \eta\, \\
& \leq &   -\frac{2}{ p}\int |\nabla (\eta f^{p/2})|^2
+\frac{2}{ p}\int |\nabla \eta|^2 f^p,
\ean
where $\nabla=\nabla_{g(t)}$.
Therefore

\ba \label{timederivative1}
\frac{\partial}{\partial t} \int f^p\eta^2
 +2
\int
   |\nabla (\eta f^{p/2})|^2
\le 2\int |\nabla \eta|^2 f^p +
(pa +\frac{\gamma}{2})\int f^p \eta^2~.
\ea

Next we define for  $0<\tau<\tau'<T$
\[
\psi(t)=\left\{ \begin{array}{ll} 0&  0\leq t\leq\tau~,\\
(t-\tau)/(\tau'-\tau)&  \tau\leq t\leq \tau'~,\\
1&  \tau'\leq t\leq T~. \end{array} \right.
\]
Multiplying (\ref{timederivative1}) by $\psi$, we obtain
\[
\frac{\partial}{\partial t}\left(\psi\int f^p\eta^2\right)+2\psi
\int|\nabla (\eta f^{p/2})|^2\leq 2\psi\int|\nabla\eta|^2 f^p
+((pa+\frac{\gamma}{2})\psi+\psi') \int f^p\eta^2~.
\]
Integrating this with respect to $t$ we get
\[
\int_t f^p\eta^2+2\int^t_{\tau'}\int|\nabla (\eta f^{p/2})|^2
\leq 2\int^T_\tau \int |\nabla\eta|^2 f^p
+\left(pa+\frac{\gamma}{2}+\frac{1}{ \tau'-\tau}\right)
\int^T_\tau \int f^p \eta^2~
\]
for $\tau' \leq t \leq T$.
Applying this estimate and the Sobolev inequality we deduce
\ba  \label{timeintegral1}
\int^T_{\tau'}\int f^{p(1+\frac{2}{ n})}\eta^{2+\frac{1}{ n}}
& \leq & \int^T_{\tau'}
\left(\int f^p\eta^2\right)^{2/n}\left(\int f^{\frac{pn}{ n-2}}
\eta^{\frac{2n}{n-2}}\right)^{\frac{n-2}{n}} \nonumber\\
&\leq &  C_{S,2}^2\Bigl(\sup_{\tau'\le t\le T} \int f^p\eta^2\Bigr)^{2/n}
\int^T_{\tau'} \int \left|\nabla(\eta f^{p/2})\right|^2  \nonumber\\
&\leq &  C_{S,2}^2\left[ 2\int^T_\tau \int
|\nabla\eta|^2 f^p+\left(pa+\frac{\gamma}{2}+\frac{1}{\tau'-\tau}\right)
\int^T_\tau \int f^p \eta^2\right]^{1+\frac{2}{n}}~.\nonumber
\\
\ea

We put
\[
H(p,\tau,R)=\int^T_\tau\int_{B(x_0, R, 0)} f^p~
\]
for $0<\tau<T$ and $0<R<d_{g(0)}(x_0, \partial M)$.
Given $0<R'<R<d_{g(0)}(x_0, \partial M)$, we
define  $\eta(x)=1$ for $d(x_0, x, 0)\ leq R'$,
$\eta(x)=1-\frac{1}{R-R'} (d(x_0, x, 0)-R')$ for
$R'\le d(x_0, x, 0)\le R$, and $\eta(x)=0$ for
$d(x_0, x, 0) \ge R$. Noticing
$|\nabla \eta|_t\le \frac{1}{R-R'}e^{-{1\over 2}\lambda_*t}$
we derive from (\ref{timeintegral1})
\be  \label{Hintegral}
H\left(p\left(1+{2\over n}\right),\tau',R'\right) \leq C_S^2
\left[pa+\frac{\gamma}{2}+\frac{1}{ \tau'-\tau}+
{2 e^{-\lambda_*T}\over (R-R')^2}\right]^{1+\frac 2n}
H(p,\tau,R)^{1+\frac 2n}~.
\ee
Now we fix $0<R<d_{g(0)}(x_0, \partial M)$ and  set
$\mu=1+{2\over n},\; p_k=p_0\mu^k,\; \tau_k=(1- {1\over \mu^{k+1}})t$
and $R_k={R\over 2}(1+{1\over \mu^{k}})$ with $R={1\over 2} dist_{g(0)}
(x,\partial N)$. Then it follows from (\ref{Hintegral}) that
\ban
\lefteqn{H(p_{k+1},\tau_{k+1}, R_{k+1})^{\frac{1}{p_{k+1}}}\leq} \\
\lefteqn{C_S^{\frac{2}{p_{k+1}}}\left[a p_{k}+\frac{\gamma}{2}
+{\mu^2\over \mu-1}\cdot {1\over
t}+
\frac{2 e^{-\lambda_*T}\mu^2}{R^2(\mu-1)^2}\right]^{\frac{1}{p_{k}}}
\mu^{\frac{k}{p_k}}H(p_k,\tau_k,R_k)^{\frac{1}{p_k}}~ \leq} \\
& & C_S^{\frac{2}{p_{k+1}}}
\left[ap_0+\frac{\gamma}{2} +\frac{\mu^2}{\mu-1}\cdot {1\over
t}+
\frac{2 e^{-\lambda_*T}\mu^2}{ R^2(\mu-1)^2}\right]^{\frac{1}{p_k}}
\mu^{\frac{2k}{p_k}}H(p_k,\tau_k,R_k)^{\frac{1}{p_k}}~.
\ean
Hence
\ban
&&H(p_{m+1},\tau_{m+1}, R_{m+1})^{\frac{1}{p_{m+1}}}\leq   C_S^{\sum^m_{0}{2\over p_{k+1}}}
\mu^{\sum_0^m\frac{2k}{p_k}} \\
&\cdot& \left[ap_0+\frac{\gamma}{2} +\frac{\mu^2}{\mu-1}\cdot {1\over
t}+
\frac{e^{-\lambda_*T}\mu^2}{R^2(\mu-1)^2}\right]^{\sum_0^m\frac{1}{p_k}}
H(p_0, \tau_0, R_0)^{\frac{1}{p_0}}.
\ean
Letting $m\to\infty$
we arrive at (\ref{linearestimate1}).
\qed
\\

\sect{Proof of Theorem A}

\noindent {\bf Proof of the estimate (\ref{RmboundA-1})} \\

 By rescaling, we can assume $r_0=1$.
Assume that the estimate (\ref{RmboundA-1}) does not hold. Then we can find for each
$1>\epsilon>0$
a Ricci flow solution $g=g(t)$ on $M\times [0, T]$ for some $M$ and $T>0$
with the properties as postulated in the statement of the theorem, such that $|Rm|(x, t)>\alpha_n t^{-1}+\epsilon^{-2}$ for some
$(x, t) \in M \times [0, T]$ satisfying $d(x_0, x, t) <\epsilon$.

We denote by $M_{\alpha_n}$ the set of pairs $(x, t)$ such that
$|Rm|(x, t)\ge \alpha_n t^{-1}$. Consider an arbitary  positive number $A>1$ such that
$(2A+1)\epsilon \le \frac{1}{2}$.
Following [Proof of Theorem 10.1, P], we choose  $(\bar x, \bar t) \in M_{\alpha_n}$ with
$0<\bar t \le \epsilon^2, d(x_0, \bar x, \bar t)<(2A+1)\epsilon$,
such that $|Rm|(\bar x, \bar t) > \alpha_n \bar t^{-1}+\epsilon^{-2}$
and
\ba \label{4Rm}
|Rm|(x, t)\le 4|Rm|(\bar x, \bar t)
\ea
whenever
\ba \label{4Rmcondition}
(x, t) \in M_{\alpha_n}, 0 <t \le \bar t, d(x_0, x, t)\le d(x_0, \bar x, \bar t)
+A|Rm|(\bar x, \bar t)^{-\frac{1}{2}}.
\ea

For the convenience of the reader, we reproduce here the argument in [Proof of Theorem 10.1, P]
for the existence of $(\bar x,\bar t)$. Let $(x_1, t_1)$ be an arbitary point in $M \times [0, T)$
such that $0<t<\epsilon^2, d(x_0, x_1, t_1) <\epsilon$ and $|Rm|(x_1, t_1) >
\alpha_n t_1^{-1}+\epsilon^{-2}$. Now if $(x_k, t_k)$ is already contructed, but  
does not have all the desired properties of $(\bar x, \bar t)$, then we can find 
$(x_{k+1}, t_{k+1})$ satisfying 
\ba \label{4kRmcondition}
(x_{k+1}, t_{k+1}) \in M_{\alpha_n}, 0 <t_{k+1} \le t_k, d(x_0, x_{k+1}, t_{k+1})&\le& d(x_0, x_k, t_k)
\nonumber \\
&+&A|Rm|(x_k, t_k)^{-\frac{1}{2}},
\ea
such that $|Rm|(x_{k+1}, t_{k+1})
>4 |Rm|(x_k, t_k)$.  It follows that $|Rm|(x_{k}, t_k) \ge 4^{k-1}|Rm|($
$x_1, t_1)\ge 
4^{k-1} \epsilon^{-2}$. By (\ref{4kRmcondition}), we then also have
$d(x_0, x_k, t_k) \le (2A+1)\epsilon.$    Since 
$g$ is smooth on $M \times [0, T)$, the former  estimate implies that the sequence $(x_k, t_k)$ must 
be finite. We can then take the last term in the sequence to be $(\bar x, \bar t)$. \\

We set $Q=|Rm|(\bar x, \bar t)$. Note that $Q>1$ because $\epsilon<1$.  \\

\noindent {\bf Claim 1} {\it  If
\ba \label{new4Rmcondition}
\bar t -\frac{1}{2}\alpha_n Q^{-1} \le t \le \bar t,
d(\bar x, x, \bar t) \le \frac{1}{10} AQ^{-\frac{1}{2}},
\ea
then}
\ba \label{distance1}
d(x_0,x,t) \le d(x_0,\bar x,\bar t)+\frac{1}{2} AQ^{-\frac{1}{2}}.
\ea
\\

Note that (\ref{distance1}) implies
\ba \label{distance2}
d(x_0, x, t) \le (2A+1)\epsilon+\frac{1}{2} AQ^{-\frac{1}{2}}\le (\frac{5}{2}A+1)\epsilon
\ea
for $(x, t)$ satisfying (\ref{new4Rmcondition}). \\

\noindent {\it Proof of Claim 1} \\

Since $(\bar x, \bar t) \in M_{\alpha_n}$, we have $Q\ge \alpha_n \bar t^{-1}$,
so
\ba \label{tbound}
\bar t- \frac{1}{2}\alpha_n Q^{-1}\ge \frac{1}{2}\bar t.
\ea

 Consider $x=\hat x$ and $t=\hat t $ satisfying (\ref{new4Rmcondition}).
  By the triangular inequality, we have
$d(x_0, \hat x, \bar t) \le d(x_0, \bar x, \bar t)+\frac{1}{10}AQ^{-\frac{1}{2}}$.
We estimate $d(x_0, \hat x, \hat t)$. For this purpose, consider
the set $I$ of $t^* \in [\hat t, \bar t]$ such that
\ba \label{1/2}
d(x_0, \hat x, t) \le d(x_0, \bar x, \bar t)+\frac{1}{2}AQ^{-\frac{1}{2}}
\ea
 for
all $t \in [t^*, \bar t]$. Obviously, $I$ is closed and $\bar t \in I$.  We claim that it is open in $[\hat t, \bar t]$.
Consider $t^* \in I$.  For each $t \in [t^*, \bar t]$,
we apply [Lemma 8.3(b),P] to $x_0, \hat x$. We set
$R=\frac{1}{2}AQ^{-\frac{1}{2}}$. For $x \in B(x_0, R, t)$
we have
$|Rm|(x, t) \le 4Q$ if $(x, t) \in M_{\alpha_n}$. If $(x, t) \not  \in
M_{\alpha_n}$, we have by (\ref{tbound})
\ba \label{2Q}
|Rm|(x, t) \le \alpha_n t^{-1} \le  2\alpha \bar t^{-1}
\le 2Q.
\ea
For $x \in
B(\hat x, R, t)$, we have $d(x_0, x, t) \le d(x_0, \hat x ,t)+d(\hat x, x, t)
\le d(x_0, \bar x, \bar t)+AQ^{-\frac{1}{2}}$. Hence $|Rm|(x, t)
\le 4Q$, if $(x, t) \in M_{\alpha_n}$.  If $(x, t) \not \in M_{\alpha_n}$,
we again obtain (\ref{2Q}).
 By [Lemma 8.3(b), P],
we have
$$
\frac{d}{dt}d(x_0, \hat x, t)\ge -2(n-1)(\frac{2}{3} \cdot 4Q \cdot \frac{1}{2}AQ^{-\frac{1}{2}}
+2A^{-1}Q^{\frac{1}{2}}) \ge -4(n-1)(A+\frac{1}{A})Q^{\frac{1}{2}}.
$$
Hence
\ban
d(x_0, \hat x, t^*) &\le& d(x_0, \hat x, \bar t) + \frac{1}{2} \alpha_n Q^{-1} \cdot 4(n-1)(A+\frac{1}{A})Q^{\frac{1}{2}} \nonumber \\ &=& d(x_0, \hat x, \bar t)+ 2(n-1)\alpha(1+\frac{1}{A^2})
AQ^{-\frac{1}{2}}) \nonumber \\   &\le& d(x_0, \bar x, \bar t)+\frac{1}{3}AQ^{-\frac{1}{2}}.
\ean
By the continuity of the distance function, the inequality (\ref{1/2})
holds true in an open neighborhood of $t^*$ in $[\hat t,
\bar t]$. It follows that $I$ is open in $[\hat t, \bar t]$. Henc we conclude
that $I=[\hat t,
\bar t]$. Consequently, we have $d(x_0, \hat x, \hat t) \le d(x_0, \bar x, \bar t)
+\frac{1}{2}AQ^{-\frac{1}{2}}$. \\

\noindent {\bf Claim 2} {\it If $(x, t)$ satisfies (\ref{new4Rmcondition}), then
the estimate (\ref{4Rm}) holds. } \\

Indeed, consider $(x, t)$ satisfying (\ref{new4Rmcondition}). If $(x, t) \in
M_{\alpha_n}$, then (\ref{distance1}) implies that the estimate (\ref{4Rm}) holds.
If $(x, t) \not \in M_{\alpha_n}$, then we have $|Rm|(x, t) \le 2Q$ as in
(\ref{2Q}). So (\ref{4Rm}) also holds. \\

Now we take $\epsilon=\frac{1}{42}$ and $A=10$. Then $\frac{1}{10}A<1$ and
$(\frac{5}{2}A+1)\epsilon=1$. So (\ref{distance2}) implies
\ba
B(\bar x, Q^{-\frac{1}{2}}, \bar t) \subset B(x_0, 1, t)
\ea
for $t \in [\bar t-\frac{1}{2}\alpha_n Q^{-1}, \bar t]$, and hence
\ba
\int_{B(\bar x, Q^{-\frac{1}{2}}, \bar t)} |Rm|^{\frac{n}{2}}(\cdot, t) dq(\cdot, t) \le \delta_0
\ea
for $t \in  [\bar t-\frac{1}{2}\alpha_n Q^{-1}, \bar t]$.  Moreover, Claim 2 implies
that the estimate (\ref{4Rm}) holds on $B(\bar x, Q^{-\frac{1}{2}}, \bar t)
\times [\bar t-\frac{1}{2}\alpha_n Q^{-1}, \bar t]$.
 We shift $\bar t$ to the time origin
and rescale $g$ by the factor $Q$ to obtain a Ricci flow solution $\bar g(t)=Qg(\bar t+Q^{-1}t)$
on $M \times [-\frac{1}{2}\alpha_n, 0]$. Now we'll deal with $\bar g$ and all the quantities will be 
associated with $\bar g$. 
We have for $\bar g$ 
\ba \label{point1}
|Rm|(\bar x, 0)=1,
\ea
and
\ba \label{new4Rm1}
|Rm|(x, t) \le 4
\ea
whenever
$$
-\frac{1}{2}\alpha_n \le t \le 0,
d(\bar x, x, 0) \le 1.
$$
Moreover there holds
\ba \label{small-A-1}
\int_{B(\bar x, 1, 0)}|Rm|^{\frac{n}{2}}(\cdot, t) dq(\cdot, t) \le \delta_0
\ea
for $t \in [-\frac{1}{2}\alpha_n,0]$.  (Note that the geodesic balls are associated with 
$\bar g$. )  By the $\kappa$-noncollapsedness
assumption, $\bar g$ is $\kappa$-noncollapsed on the scale $Q^{\frac{1}{2}}\rho \ge \rho \ge r_0=1$. 
Finally, one readily verifies that 
\ba \label{newdistance} 
d_{\bar g(0)}(\bar x, \partial M) > \frac{1}{2}.
\ea

Applying the $\kappa$-noncollapsedness property of $\bar g$ and (\ref{new4Rm1}) we deduce
\ba \label{volumelowerbound}
V(\bar x, r, 0) \ge  \kappa r^n
\ea
for all  $0<r \le \frac{1}{2}$. 
By (\ref{new4Rm1}), (\ref{volumelowerbound}), (\ref{newdistance}) and (\ref{sobolev-2}) we can apply  Theorem B.1 in Appendix B to
 infer 
\ba
C_{S, 2, \bar g (0)}(B(\bar x, \rho(n, \kappa), 0))\le C_1(n) \equiv \frac{2(n-1)}{n-2} \frac{\omega_n^{\frac{n-1}{n}}}{2^{\frac{n-1}{n}}\omega_{n-1}}
\ea
for a positive constant $\rho(n, \kappa)$ depending only on $n$ and $\kappa$. 
By the curvature bound (\ref{new4Rm1}) and the argument in [Ye1] for
the evolution of the Sobolev constant, we then infer
\ba \label{sobolevA}
C_{S,2, \bar g(t)}(B(\bar x, \rho(n, \kappa),0)) \le C_2(n)
\ea
for $t \in [-\frac{1}{2}\alpha_n, 0]$, where $C_2(n)$
is a positive constant depending only $n$.

On the other hand, it is easy to see that $\rho(n, \kappa) \le \frac{1}{16}$. Hence the curvature bound (\ref{new4Rm1}) and
the Ricci flow equation imply that $B(\bar x, \rho(n, \kappa), 0)
\subset B(\bar x, 1, t)$ for $t \in [-\bar \alpha_n, 0]$, where
$\bar \alpha_n\le \frac{1}{2}\alpha_n$ is a positive constant depending only on $n$.
It follows that
\ba \label{smallA}
\int_{B(\bar x, \rho(n, \kappa), 0)} |Rm|^{\frac{n}{2}}(\cdot, t)dq(\cdot, t) \le
\delta_0
\ea
for $t\in [-\bar \alpha_n, 0]$.

Now we appeal to the evolution equation of $Rm$ associated with the Ricci flow
\ba \label{curvatureevolution}
\frac{\partial Rm}{\partial t}=\Delta Rm +B(Rm,Rm),
\ea
where $B$ is a certain quadratic form. It implies 
\ba \label{curvaturenorm}
\frac{\partial}{\partial t}|Rm| \le \Delta |Rm| + c(n) |Rm|^2
\ea
for a positive constant $c(n)$ depending only on $n$.
 On account of (\ref{new4Rm1}), (\ref{sobolevA}) and
(\ref{smallA}) we can apply Theorem 2.1 to (\ref{curvaturenorm}) with $p_0=\frac{n}{2}$  to deduce
\ba \label{RmfinalA-1}
|Rm|(\bar x, 0)
&\leq& (1+\frac{2}{n})^{\frac{2\sigma_n}{n}} C_2(n)^2 C_3(n)
\Bigl(\int^0_{-\bar \alpha_n}\int\limits_{B(x_0, \rho(n, \kappa),0)} |Rm|^{\frac{n}{2}}(\cdot, t)dq(\cdot, t)
\Bigr)^{\frac{2}{n}} \nonumber \\
&\le& (1+\frac{2}{n})^{\frac{2\sigma_n}{n}} C_2(n)^2 C_3(n, \kappa)
(\bar \alpha_n \delta_0
)^{\frac{2}{n}},
\ea
where
\ban
C_3(n, \kappa)=\Bigl(2c(n)n+4n(n-1) +\frac{n}{2}(1+\frac{n}{2})^2 \cdot \frac{1}{\bar \alpha_n}
+
\frac{(n+2)^2}{4\rho(n, \kappa)^2}e^{8(n-1) \bar \alpha_n}\Bigr)^{n+2\over n}.
\ean
We deduce $|Rm|(\bar x, 0) \le \frac{1}{2}$, provided that we define
\ban
\delta_0 &=& \frac{1}{2^{\frac{n}{2}}}(1+\frac{2}{n})^{-\sigma_n} C_2(n)^{-n} C_3(n, \kappa)^{-\frac{n}{2}}
\bar \alpha_n^{-1}.
\ean
But this contradicts (\ref{point1}). Hence the estimate (\ref{RmboundA-1}) has been
established. \\

\noindent {\bf Proof of the estimate (\ref{RmboundintA-1}) } \\

Consider a fixed $t_0 \in (0, T)$. If  the ratio $\frac{r_0^2}{t_0} \ge 1$,
we rescale $g$ by $t_0^{-1}$. If $\frac{r_0^2}{t_0} \le 1$, we rescale
$g$ by $r_0^{-2}$. We handle the former case, while the latter is similar.
For the rescaled flow $t_0^{-1}g(t_0t)$ on $[0, t_0^{-1}T)$ we have by
(\ref{RmboundA-1})
\ba
|Rm|(x, t) \le 2\alpha_n+\epsilon_0^{-2}
\ea
whenever $\frac{1}{2} \le t <t_0^{-1}T$ and $d(x_0, x, t) <
t_0^{-\frac{1}{2}}r_0 \epsilon_0$. We rescale the flow further by the factor $\lambda_n\equiv 2\alpha_n+\epsilon_0^{-2}
>42$
to obtain $\bar g(t)=\lambda_n t_0^{-1} g(\lambda_n^{-1}t_0t)$ on $[0, \lambda_n t_0^{-1}T)$.
Note that the time $t_0$ is transformed to the time 
$\lambda_n$.  We deal with $\bar g$ and all quantities will be associated with $\bar g$.  We have for $\bar g$
\ba
|Rm|(x, t) \le 1
\ea
whenever $\frac{\lambda_n}{2} \le t < \lambda_n t_0^{-1}T$ and
$d(x_0, x, t) <\sqrt{\lambda_n} t_0^{-\frac{1}{2}} r_0 \epsilon_0 $. Note
$\sqrt{\lambda_n} t_0^{-\frac{1}{2}} r_0 \epsilon_0>1$ and $\sqrt{\lambda_n} t_0^{-\frac{1}{2}}\rho 
\ge \sqrt{\lambda_n} t_0^{-\frac{1}{2}} r_0 \ge \sqrt{\lambda_n}>6$.  It follows in particular 
that $\bar g$ is $\kappa$-noncollapsed on the scale $6$. We also have 
$d_{\bar g(t)}(x_0, \partial M) \ge \lambda_n$ for all $t \in [0, \lambda_n t_0^{-1}T)$. 
We deduce
\ba
V(x, r, \lambda_n) \ge \kappa r^n
\ea
whenever 
$d(x_0, x, \lambda_n) \le \frac{1}{2}$ and $0<r \le \frac{1}{2}$.    
Now we can argue as
in the above proof of the estimate (\ref{RmboundA-1}) to infer
$B(x, \frac{1}{4}, \lambda_n) \subset B(x, \frac{1}{2}, t)$ and  
\ba
C_{S, 2, \bar g(t)}(B(x, \rho(n, \kappa), \lambda_n)) 
\le  C_2(n)
\ea
whenever $\lambda_n-\bar \alpha_n \le t \le  \lambda_n$ 
and $d(x_0, x, \lambda_n) \le \frac{1}{2}$, where 
$\bar \alpha_n>0$ and $C_2(n)>0$ depend only on $n$, and  
$0<\rho(n, \kappa) \le \frac{1}{16}$ depends  only on $n$ and $\kappa$.  
(We retain the notations $\bar \alpha_n, \rho(n, \kappa)$ and $C_2(n)$ although they may be different 
from their values before.) 
Now we can  apply Theorem 2.1 on
$B(x, \rho(n, \kappa), \lambda_n) \times [\lambda_n-\bar \alpha_n, \lambda_n]$ for 
each $x \in B(x_0, \frac{1}{2}, \lambda_n)$ to 
derive (for $\bar g$)
\ba
|Rm|(x, \lambda_n) &\le& \tilde C_0(n, \kappa)
\left(\int_{\lambda_n-\bar \alpha_n}^{\lambda_n} \int_{B(x, \rho(n, \kappa), \lambda_n)}|Rm|^{\frac{n}{2}}(\cdot, t)dq(\cdot, t)dt\right)^{\frac{2}{n}} \nonumber \\
&\le&  C_0(n, \kappa) 
\sup\limits_{\lambda_n-\bar \alpha_n \le t
\le \lambda_n} \left(\int_{B(x, \frac{1}{2}, t)}|Rm|^{\frac{n}{2}}(\cdot, t)dq(\cdot, t)\right)^{\frac{2}{n}}.
\ea
for all $x\in B(x_0, \frac{1}{2}, \lambda_n)$, where $\tilde C_0(n, \kappa)>0$ depends only 
on $n$ and $\kappa$ and $C_0(n, \kappa)=\bar \alpha_n \tilde C_0(n, \kappa)$.  
Scaling back to $g$ we then arrive at the desired estimate (\ref{RmboundintA-1})
(with $t_0$ in place of $t$). \\
\qed \\

\sect{Proof of Theorem B}

\noindent {\bf Proof of the estimate (\ref{RmboundB-1})} \\

Assume that the estimate (\ref{RmboundB-1}) fails to hold.
Then we carry out the same construction as in the proof of Theorem A.
Again we assmue $r_0=1$ and 
choose $\epsilon=\frac{1}{42}$ and $A=10$. We deal with the rescaled flow $\bar g$ and 
all quantities will be associated with $\bar g$.
By (\ref{smallRmintegralB-1}) we have for $\bar g$, in place of (\ref{small-A-1}) 
\ba \label{small-B-1}
\int_{B(\bar x, r, t)} |Rm|^{\frac{n}{2}}(\cdot, t)dq(\cdot, t)\le r^{-n} \delta_0  V(\bar x, r, t)
\ea
for all $0<r\le \frac{1}{2}$ and $t \in (-\frac{1}{2}\alpha_n, 0].$  As before, we also have for $\bar g$
\ba \label{pointB-1}
|Rm|(\bar x, 0)=1
\ea
and
\ba \label{4RmB-1}
|Rm|(x, t) \le 4
\ea
whenever
\ba
-\frac{1}{2}\alpha_n\le t \le 0, d(\bar x, x, 0) \le 1.
\ea
Moreover, we have 
\ba \label{newdistance1} 
d_{\bar g(0)}(\bar x, \partial M) > \frac{1}{2}
\ea
and 
\ba \label{newdiameter}
diam_{\bar g(t)}(M) \ge 1
\ea
for all $ t \in [-\frac{1}{2}\alpha_n, 0]$. By (\ref{4RmB-1}), (\ref{newdistance1})  and the Ricci flow equation 
we have
\ba \label{newdistance2}
d_{\bar g(t)}(\bar x, \partial M) \ge \frac{1}{3}
\ea
for all $t \in [-\alpha_n^*, 0]$, where $\alpha_n^* \le \frac{1}{2}\alpha_n$ depends 
only on $n$. 

By (\ref{4RmB-1}), (\ref{newdiameter}) and (\ref{newdistance2}) we can apply Theorem B.2 in Appendix B
 to deduce
\ba \label{sobolevB-1}
C_{S, 2, \bar g(t)}(B(\bar x, \frac{1}{12}, t)) \le \frac{C_5(n)}{V(\bar x, \frac{1}{12}, t)^{\frac{1}{n}}}
\ea
for $t \in [-\alpha_n^*, 0]$, with a positive constant $C_5(n)$
depending only on $n$.  On the other hand, (\ref{4RmB-1}) implies that
\ba \label{includeB-1}
B(\bar x, \frac{1}{14}, t_1) \subset  B(\bar x, \frac{1}{13}, t_2) \subset
B(\bar x, \frac{1}{12}, t_3)
\ea
for all $t_1, t_2$ and $t_3 \in [-\bar \alpha_n, 0]$, with a positive constant $\bar \alpha_n
\le \min\{\alpha_n^*, \frac{1}{2} \alpha_n\}$ depending only on $n$. Consequently, we have
\ba \label{smallintegralB-1-2}
\int_{B(\bar x, \frac{1}{14}, 0)} |Rm|^{\frac{n}{2}}(\cdot, t)dq(\cdot, t) \le (13)^n \delta_0  vol_{\bar g(t)}(B(\bar x, \frac{1}{13}, t))
\ea
and
\ba \label{sobolevB-1-2}
C_{S, 2, \bar g(t)}(B(\bar x, \frac{1}{14}, 0)) \le \frac{C_5(n)}{V(\bar x, \frac{1}{12}, t)^{\frac{1}{n}}}
\ea
for all $t\in [-\bar \alpha_n, 0]$. Moreover, (\ref{includeB-1}) combined
with (\ref{4RmB-1}) leads via the Ricci flow equation to
\ba \label{volumeB-1-2}
\min\limits_{-\bar \alpha_n \le t \le 0} V(\bar x, \frac{1}{12}, t)
\ge e^{-4n(n-1)\bar \alpha_n} \max\limits_{-\bar \alpha_n\le t \le 0} V(\bar x, \frac{1}{13}, t)
\ea
for each $t \in [-\bar \alpha_n, 0]$.
Now we apply Theorem 2.1 to deduce
\ba \label{RmfinalB-1}
|Rm|(\bar x, 0)
&\leq& (1+\frac{2}{n})^{\frac{2\sigma_n}{n}} \frac{C_5(n)^2C_6(n)}{\min\limits_{-\bar \alpha_n \le t \le 0} V(\bar x, \frac{1}{12}, t)^{\frac{2}{n}}} 
\Bigl(\int^0_{-\bar \alpha_n}\int\limits_{B(\bar x, \frac{1}{14},0)} |Rm|^{\frac{n}{2}}(\cdot, t)dq(\cdot, t)
\Bigr)^{\frac{2}{n}} \nonumber \\
&\le& \frac{169(1+\frac{2}{n})^{\frac{2\sigma_n}{n}}
C_5(n)^2C_6(n)(\bar \alpha_n \delta_0)^{\frac{n}{2}}}{\min\limits_{-\bar \alpha_n \le t \le 0} V(\bar x, \frac{1}{12}, t)^{\frac{2}{n}}}  \max\limits_{-\bar \alpha_n
\le t \le 0} V(\bar x, \frac{1}{13}, t)^{\frac{2}{n}}
 \nonumber \\
&\le& 169(1+\frac{2}{n})^{\frac{2\sigma_n}{n}} C_5(n)^2 C_6(n)
(\bar \alpha_n \delta_0
)^{\frac{2}{n}}e^{8(n-1)\bar \alpha_n},
\ea
with a positive constant $C_6(n)$ depending only on $n$. Choosing
\ban
\delta_n=\frac{1}{2^{\frac{n}{2}}(13)^n}(1+\frac{2}{n})^{-\sigma_n} C_5(n)^{-n} C_6(n)^{-\frac{n}{2}}
\bar \alpha_n^{-1} e^{-8(n-1)\bar \alpha_n}
\ean
we then obtain $|Rm|(\bar x, 0)\le \frac{1}{2}$, contradicting
(\ref{pointB-1}).
\\

\noindent {\bf Proof of the estimate (\ref{RmboundintB-1})} \\

Consider a fixed $t_0 \in (0, T)$. As in the corresponding part of
the proof of Theorem A, we present the case $\frac{r_0^2}{t_0} \ge 1$,
and rescale $g$ by $t_0^{-1}$. Again we have for
the rescaled flow $t_0^{-1}g$ on $[0, t_0^{-1}T)$
\ba
|Rm|(x, t) \le 2\alpha_n+\epsilon_0^{-2}
\ea
whenever $\frac{1}{2} \le t <t_0^{-1}T$ and $d(x_0, x, t) <
t_0^{-\frac{1}{2}}r_0 \epsilon_0$. As before, we rescale the flow further by the factor $\lambda_n\equiv 2\alpha_n+\epsilon_0^{-2}>42$
to obtain $\bar g(t)=\lambda_n t_0^{-1} g(\lambda_n^{-1}t_0t)$ on $[0, \lambda_n t_0^{-1}T)$.
Again, the time $t_0$ is transformed to the time 
$\lambda_n$. As before we deal with $\bar g$ and all quantotoes are associated with 
$\bar g$. We have $\sqrt{\lambda_n} 
t_0^{-\frac{1}{2}}r_0\epsilon_0>1$ and the curvature estimate for $\bar g$
\ba \label{Rm-one}
|Rm|(x, t) \le 1
\ea
whenever $\frac{\lambda_n}{2} \le t < \lambda_n t_0^{-1}T$ and
$d(x_0, x, t) \le \sqrt{\lambda_n} t_0^{-\frac{1}{2}} r_0 \epsilon_0$.
We also have $d_{\bar g(t)}(x_0, \partial M) \ge \lambda_n$ and 
$diam_{\bar g(t)} (M) \ge \lambda_n $  for 
all $t \in [0, \lambda_n t_0^{-1}T)$. 

Now we can apply the arguments in  
the above proof of (\ref{RmboundB-1}).  Since the conition $d_{\bar g(t)}(x_0, \partial M) \ge 
\frac{1}{3}$ is now replaced by $d_{\bar g(t)}(x_0, \partial M) \ge \lambda_n >1$ and we have 
$\sqrt{\lambda_n} t_0^{-\frac{1}{2}} r_0 \epsilon_0>1$ (for the purpose of applying 
(\ref{Rm-one})), we can 
replace the radii $\frac{1}{14}, \frac{1}{13}$ and $\frac{1}{12}$ by $\frac{1}{10}, 
\frac{1}{9}$ and $\frac{1}{8}$.  We deduce
\ba
|Rm|(x, \lambda_n) &\le& \frac{C_7(n)}{\min \limits_{\lambda_n-\bar \alpha_n \le t \le \lambda_n} V(x, \frac{1}{8}, t)^{\frac{2}{n}}} \left( \int_{\lambda_n-\bar \alpha_n}^{\lambda_n} \int_{B(x, \frac{1}{10}, \lambda_n)} |Rm|^{\frac{n}{2}}(\cdot, t)dq(\cdot, t)\right)
^{\frac{2}{n}}  \nonumber \\
&\le & C_0(n)
\sup\limits_{\lambda_n-\bar \alpha_n \le t
\le \lambda_n} \left(\frac{
\int_{B(x, \frac{1}{9}, t)}|Rm|^{\frac{n}{2}}(\cdot, t)dq(\cdot, t)}{V(x, \frac{1}{9}, t)}\right)^{\frac{2}{n}},
\ea
 with positive constant constants 
$C_7(n), \bar \alpha_n$ and $C_0(n)$ depending only on $n$, whenever $d(x_0, x, \lambda_n)
\le \frac{1}{2}$.  
Scaling back to $g$ we then arrive at the desired estimate (\ref{RmboundintB-1})
(with $t_0$ in place of $t$).  \\

\sect{Proof of Theorem C}

\noindent {\bf Proof of Theorem C}
\\

We establish the condition (\ref{smallRmintegralB-1}). Then the theorem follows from
Theorem B.
By rescaling we can assume $r_0=1$. Then (\ref{riccibound}) becomes
\ba \label{ricciboundC-1}
Ric \ge -(n-1)g.
\ea
By (\ref{smallRmintegralC-1}), we have now
\ba \label{smallintegralC-1}
\int_{B(x,1,t)} |Rm|^{\frac{n}{2}}(\cdot, t)dq(\cdot, t) \le \delta_0 V(x,1,t)
\ea
for all $t \in [0, T)$.
 By Bishop-Gromov relative volume comparison, we have
\ba \label{volumeC-1}
V(x, R, t) \le \frac{V_{-1}(R)}{V_{-1}(r)} V(x, r, t)
\le C(n) \frac{V(x, r, t)}{r^n},
\ea
with a positive constant $C(n)$ depending only on $n$, provided that $t \in [0, T],
d(x_0,x, t)<1$, and $0<r<R\le 1-d(x_0,x,t)$. Here $v_{-1}(r)$ denotes the volume
of a geodesic ball of radius $r$ in $\mbox{\bf H}^n$, the $n$-dimensional
hyperbolic space
(of sectional curvature $-1$). If $t \in [0, T)$ and $ d(x_0, x,t) \le \frac{1}{4}$,
we then have $B(x_0, \frac{1}{4}, t) \subset B(x, \frac{1}{2}, t)
\subset B(x_0, 1, t)$. Consequently,
\ba
V(x, r,t) &\ge& C(n)^{-1} r^n V(x, \frac{1}{2}, t)
\ge C(n)^{-1} r^n V(x_0, \frac{1}{4}, t) \nonumber \\  &\ge&
4^{-n} C(n)^{-2} r^n V(x_0, 1, t)
\ea
for $0<r\le \frac{1}{2}$. This leads to 
\ba
V(x, r, t) \ge 4^{-n}C(n)^{-2} \frac{r^n}{r_0^n} V(x_0, r_0, t)
\ea
as long as $t \in [0, T), d(x_0, x, t) \le \frac{1}{4}r_0$ and 
$0<r\le \frac{1}{2}r_0$. 
Hence we infer 
\ba \label{smallintegralC-1-2}
\int_{B(x, r, t)} |Rm|^{\frac{n}{2}}(\cdot, t)dq(\cdot, t)  &\le& \int_{B(x_0, r_0, t)}
|Rm|^{\frac{n}{2}}(\cdot, t) dq(\cdot, t) \le \delta_0 \frac{V(x_0, r_0, t)}{r_0^n}
\nonumber \\ &\le& 4^n C(n)^2 \delta_0 \frac{V(x, r ,t)}{r^n}
\ea
whenever $t \in [0,T), d(x_0, x, t) \le \frac{1}{4} r_0$ and $0<r\le \frac{1}{2}r_0$.
Choosing $\delta_0$ to be the $\delta_0$ in Theorem B multiplied
by $4^{-n}C(n)^{-2}$ and replacing $r_0$ by $\frac{r_0}{2}$ we
then have all the conditions of Theorem B. The desired estimate follows.\\
\qed
\\

\vspace{1cm}

\noindent  {\bf {\Large Appendices}} \\

\appendix
\sect{$\kappa$-Noncollapsedness}

In this appendix we present a stronger version of [Theorem 4.1, P], namely Theorem A.1 regarding 
the $\kappa$-noncollapsed property of the Ricci flow at finite times. We observed 
this version in the beginning of 2003. Later, we learnt that Perelman also made 
the same observation, see [KL].  \\

\noindent {\bf Definition} Let $g$ be a Riemannian metric on a manifold $M$ of dimension
$n$. Let $\kappa$ and $\rho$ be positive numbers. We say that $g$ is {\it $\kappa$-noncollapsed on the scale 
$\rho$ relative to the positive part of the scalar curvature (or relative to upper bounds of the 
scalar curvature)},
if $g$ satisfies  $vol(B(x, r))\ge \kappa r^n $
for all $x \in M$ and $r>0$ satisfying $r < \rho$ and $\sup\{R(x): x\in B(x, r)\}
\le r^{-2}$. (Note that $\sup\{R(x): x\in B(x, r)\}
\le r^{-2}$ is equivalent to $\sup\{R^+(x): x\in B(x, r)\}
\le r^{-2}$. Hence the terminology 
``the positive part of the scalar curvature".)  We say that a family of Riemannian metrics $g=g(t)$ is $\kappa$-noncollapsed on the scale 
$\rho$ relative to the positive part of the scalar curvature, if $g(t)$ is $\kappa$-noncollapsed on the scale $\rho$ 
relative to the positive part of the scalar curvature for each $t$ (in the given domain). 
\\

Obviously, if $g$ is  $\kappa$-noncollapsed on the scale 
$\rho$ relative to the positive part of the scalar curvature, then it is 
$\tilde \kappa$-noncollapsed on the scale 
$\rho$, where $\tilde \kappa=c(n)\kappa$ for a positive constant 
$c(n)$ depending only on the dimension $n$. 

\begin{theo} Let $g=g(t)$ be a smooth solution of the Ricci flow on
$M \times [0, T)$ for a closed manifold $M$ of dimension $n\ge 2$ and some
finite $T>0$. Let $\rho>0$ be an arbitary positive number. Then $g$ is $\kappa$-noncollapsed on the scale $\rho$ relative to 
the positive part of the scalar curvature
for $t\in [0,T)$, where $\kappa=\kappa(T+\rho^2, g(0))$ depends  on 
(an upper bound of) $T+\rho^2$ and the initial metric $g(0)$.  
\end{theo}
\Pf Let $g=g(t)$ be a smooth solution of the Ricci flow on
$M \times [0, T)$ for a closed manifold $M$ of dimension $n\ge 2$ and some
finite $T>0$.  Assume that there is no $\kappa>0$ such that 
$g$ is $\kappa$-noncollapsed on the scale $\rho$ relative to 
the positive part of the scalar curvature for $t \in [0,T)$. Then 
there is a sequence of times $t_k \in [0, T)$ with $t_k \rightarrow T$, a 
sequence of points $x_k \in M$, and a sequence of positive 
numbers $r_k <\rho$ such that for all $k$
\ba \label{smallvolume1}
V(x_k, r_k, t_k) \le \frac{1}{2k} \omega_n r_k^n
\ea 
and 
\ba \label{scalarbound}
R(\cdot, t_k) \le r_k^{-2}
\ea
on $B(x_k, r_k, t_k)$.  Here $\omega_n$ denotes the volume of the 
$n$-dimensional  euclidean ball  of radius 1. 

For a fixed $k$ we set $\rho_{k,j}=\frac{r_k}{2^j}$.  Choose the largest 
$j$ such that $V(x_k, \rho_j, t_k) \le \frac{1}{2k} \omega_n \rho_j^n$.
Since $\frac{V(x_k, \rho, t_k)}{\omega_n\rho^n} \rightarrow 1$ as  $\rho \rightarrow 0$,
such a $j$ exists. Then we have $V(x_k, \rho_{k, j+1}, t_k) > \frac{1}{2k} \omega_n \rho_{k,j+1}^n$.
We replace the value of $r_k$ by $\rho_{k,j}$ with this largest $j$. Then we have in addition to 
(\ref{smallvolume1}) and (\ref{scalarbound})
\ba \label{smallvolume2}
V(x_k, \frac{r_k}{2}, t_k) >  \frac{1}{2k} \omega_n (\frac{r_k}{2})^n.
\ea

Next we consider the entropy functional of Perelman [P]
\ba
W(g, f, \tau) &=& \int_M \left[ \tau(|\nabla f|^2 +R)+
f-n\right]  (4\pi\tau)^{-\frac{n}{2}} e^{-f}dvol
\ea
for smooth metrics $g$ and Lipschitz functions $f$ on $M$, and $\tau>0$, under the side condition
\ba \label{side}
(4\pi \tau)^{-\frac{n}{2}}\int_M e^{-f}dvol =1,
\ea
where all geometric quantities are associated with $g$.
We construct a sequence of Lipschitz  functions 
$f_k$ on $M$ satisfying the side condition (\ref{side}) for $g=g(t_k)$ such that $W(g(t_k), f_k, r_k^2) \rightarrow -\infty$. 
For $\delta>0$ we set 
$\psi_{\delta}(t)=1$ for $0\le t \le \frac{1}{2}$, $\psi_{\delta}(t)=\delta$ for $t \ge 1$, 
and $\psi_{\delta}(t)= \frac{2(1+\delta-t)}{1+2\delta}$ for $\frac{1}{2}\le t \le 1$. Let $\Lambda>0$.
Following [P]  we then set $f_k=-2 \log(\Lambda \psi_{\delta}(\frac{r}{r_k}))$, where $r(x)=d(x_k, x, t_k)$.   
Obviously, 
\ba
\int_{B(x_k, \frac{r_k}{2}, t_k)}f_ke^{-f_k}(4\pi r_k^2)^{-\frac{n}{2}}dvol_{g(t_k)} =
-2 (4\pi)^{-\frac{n}{2}} \frac{V(x_k, \frac{r_k}{2}, t_k)}{r_k^n} \Lambda^2 \log \Lambda.
\ea
By the conditions (\ref{smallvolume1}), (\ref{smallvolume2}) and (\ref{scalarbound}), this 
integral is  
dominating in $W(g(t_k), f_k, r_k^2)$, provided that $\delta$ is small and $\Lambda$ is large. 
It is easy to choose $\Lambda=\Lambda_k \approx \omega_n^{-1} (4\pi)^{\frac{n}{2}} k 2^{n+1} $ 
and $\delta=\delta_k$ with $\delta_k^2 \log \delta_k \approx -\frac{1}{2}r_k^n vol_{g(t_k)}(M)^{-1}$
such tht 
\ba
W(g(t_k), f_k, r_k^2) <-\log \Lambda_k
\ea
and 
\ba
(4\pi)^{-\frac{n}{2}} \int_M e^{-f_k} dvol_{g(t_k)}=1.
\ea

For a fixed $k$ let $\tau(t)=t_k-t+r_k^2$ and $f$ be the solution of the equation
\ba
\frac{\partial f}{\partial t}=-\Delta f+|\nabla f|^2-R+\frac{n}{2\tau}
\ea
on $[0, t_k]$ associated  with $g=g(t)$, with the initial value $f(\cdot, t_k)=f_k$. 
By the monotonicity of  the entropy functional [P], 
we have 
for $\bar f_k=f(\cdot, 0)$
\ba \label{newentropy}
W(g(0), \bar f_k, t_k+r_k^2) \le W(g(t_k), f_k, r_k^2) <-\log \Lambda_k. 
\ea
Moreover, we also have 
\ba
(4\pi (t_k+r_k^2))^{-\frac{n}{2}} \int_M e^{-\bar f_k} dvol_{g(0)}=1.
\ea
However, $t_k+r_k^2 \le T+\rho^2$. Hence the logarithmic Sobolev inequality implies that 
$W(g(0), \bar f_k, t_k+r_k^2)$ is bounded from below by a finite constant independent of $k$ (see 
[R] and [Y4]). This contradicts (\ref{newentropy}). \qed

\sect{Estimates of the Sobolev Constant}

In this appendix we present two  estimates for the Sobolev
constant.

\begin{lem}    Let $(M, g)$ be a Riemannian manifold of dimension $n$.
Assume that the sectional curvatures $K_g$ of $g$ satisfies 
$\kappa_1 \le K_g \le \kappa_2$ on a geodesic ball $B(p, r_0)$ in $(M, g)$, 
such that $r_0 \le d(p, \partial M)$. Set $r_1=\frac{1}{4}\min\{r_0, \frac{\pi}{4\sqrt{\kappa_2}}\}$. 
Then the injectivity radius $i(q)$ at any $q 
\in B(p, r_1)$ 
satisfies 
\be 
i(q) \ge r_2, \tag{B.1} 
\ee
where 
\be
r_2=\frac{r_1}{2} \left(1+\frac{V_{\kappa_1}(2r_1)^2}{vol_g(B(p, r_1))V_{\kappa_1}(r_1)}\right)^{-1}, \tag{B.2}
\ee
 $r_1=\frac{1}{4}\min\{r_0, \frac{\pi}{\sqrt{\kappa_2}}\}$,  and for any $r>0$, $V_{\kappa_1}(r)$ denotes the 
volume of a geodesic ball of radius $r$ in the $n$-dimensional model space (a simply connected 
complete Riemannian manifold) of sectional curvature $\kappa_1$. 
\end{lem}
\Pf The estimate (B.1) follows from the proof of a similar estimate [(4.23), CGT]  in  [CGT]. In [CGT], 
the estimate [(4.23), CGT] is proved for complete manifolds. Under the above assumption 
about the radius $r_0$, one can easily check that the involved geodesics in the relevant arguments in [CGT] 
all stay inside $B(p, r_0)$, and hence those arguments all go through. \qed

\begin{theo} Let $(M, g)$ be a Riemannian manifold of dimension $n$.
Assume that the sectional curvature $K_g$ of $g$ satisfies 
$\kappa_1 \le K_g \le \kappa_2$ on a geodesic ball $B(p, r_0)$ in $(M, g)$, 
such that $r_0 \le d(p, \partial M)$. Let $r_1$ and $r_2$ be defined as in Lemma A.1.
 Then we have 
\be
C_{S,g}(B(x, r)) \le  \frac{\omega_n^{\frac{n-1}{n}}}{2^{\frac{n-1}{n}}\omega_{n-1}}
\tag{B.3}
\ee
for all $0<r \le \frac{r_2}{2}$, 
where $\omega_n$ denotes the volume of the unit sphere in the $(n+1)$-dimensional Euclidean 
space. 
\end{theo}
\Pf This follows from Lemma B.1 and Croke's isoperimetric inequality, see  [C] or [Proposition 
V.2.3(a), Ch]. \qed

\begin{theo} Let $(M, g)$ be a Riemannian manifold of dimension $n$. Assume that 
the Ricci curvature $Ric_g$ of $g$ satisfies $Ric_g \ge -\kappa(n-1)$ for a nonnegative 
constant $\kappa$
on a geodesic ball $B(p, r_0)$ in $(M, g)$, 
such that $r_0 \le \min\{d(p, \partial M), diam(M)\}$. 
Then we have 
\be
C_{S, g}(B(p, r)) \le C(n, \kappa) \left(\frac{V_{-\kappa}(r)}{vol_g(B(p, r))}\right)^{\frac{1}{n}}
\tag{B.4}
\ee
for all $0<r\le \min\{\frac{r_0}{4}, 1\}$,
where $C(n, \kappa)$ is a positive constant depending only on $n$ and $\kappa$.
\end{theo}
\Pf The estimate (B.4) was established in [An] for complete Riemannian manifolds. 
Its proof in [An], which is of local nature, can be carried over to the present situation, because 
the involved geodesics all stay inside of $B(p, r_0)$.  Note that the proof in [An] 
is based on some inequalities established in [G]. Those inequalities also carry over 
to the present situation for the same reason. \qed \\

\noindent {\bf Remark} Obviously, Theorem B.3 leads to a different estimate of the Sobolev constant 
under the conditions of Theorem B.2. This estimate can replace Theorem B.2 in the proof of Theorem A.

\end{document}